\documentclass{amsart}

\theoremstyle{plain} \numberwithin{equation}{section}
\newtheorem{Theorem}{Theorem}
\newtheorem{Lemma}[Theorem]{Lemma}

\newtheorem{Proposition}[Theorem]{Proposition}

\theoremstyle{remark}

\title[Bari-Markus property]
{Bari-Markus property for Riesz projections of Hill operators with
singular potentials}

\author{Plamen Djakov}

\author{Boris Mityagin}

\begin{document}

\address{Sabanci University, Orhanli,
34956 Tuzla, Istanbul, Turkey}

 \email{djakov@sabanciuniv.edu}

\address{Department of Mathematics,
The Ohio State University,
 231 West 18th Ave,
Columbus, OH 43210, USA}

\email{mityagin.1@osu.edu}

\begin{abstract}
The Hill operators $L y = - y^{\prime \prime} + v(x) y, \; x \in
[0,\pi],$ with $H^{-1}$ periodic potentials, considered with
periodic, antiperiodic or Dirichlet boundary conditions, have
discrete spectrum, and therefore, for sufficiently large $N,$ the
Riesz projections
$$
 P_n = \frac{1}{2\pi i} \int_{C_n}
(z-L)^{-1} dz, \quad C_n=\{z: \; |z-n^2|= n\}
$$
are well defined. It is proved that $$\sum_{n>N} \|P_n -
P_n^0\|^2_{HS} < \infty,
$$
where $P_n^0$ are the Riesz projection of the free operator and
$\|\cdot\|_{HS}$ is the Hilbert--Schmidt norm.
\end{abstract}

\subjclass[2000]{34L40 (primary), 47B06, 47E05 (secondary)}

\maketitle

\section{Introduction}

We consider the Hill operator
\begin{equation}
\label{01} Ly = - y^{\prime \prime} + v(x) y, \qquad x \in
I=[0,\pi],
\end{equation}
with a singular complex--valued periodic potential $v, \; v(x+\pi) =
v(x), \;v \in H^{-1}_{loc} (\mathbb{R}),$ i.e., $$ v(x) = v_0 +
Q^\prime (x),$$ where $$Q \in L^2_{loc}(\mathbb{R}), \quad Q(x+\pi)
= Q(x), \quad q(0) = \int_0^\pi Q(x) dx =0,$$ so
\begin{equation}
\label{001} Q = \sum_{m \in 2\mathbb{Z}\setminus \{0\}} q (m) e^{
imx},\quad \|v|H^{-1} \|^2 = |v_0|^2 + \sum_{m \in
2\mathbb{Z}\setminus \{0\}} |q(m)|^2/m^2 <\infty.
\end{equation}

A. Savchuk and A. Shkalikov \cite{SS03} gave thorough spectral
analysis of such operators.  In particular, they consider a broad
class of boundary conditions (bc) -- see (1.6), Theorem 1.5 there --
in terms of a function $y$ and its quasi--derivative
$$
u = y^\prime - Q y.
$$
The natural form of periodic or antiperiodic $(Per^\pm)$ bc is the
following one:
\begin{equation}
\label{02}
 Per^\pm: \quad y(\pi) = \pm y(0), \quad u(\pi) = \pm u(0)
\end{equation}
If the potential $v$ happens to be an $L^2$-function these $bc$ are
identical to the classical ones (see discussion in \cite{DM16},
Section 6.2).

 The Dirichlet bc is more simple:
$$ Dir: \quad  y(0) =0,  \quad y(\pi) =0; $$
it does not require quasi--derivatives, so it is defined in the same
way as for $L^2$--potentials $v$.

In our analysis of instability zones of Hill and Dirac operators
(see \cite{DM15} and the comments there) we follow an approach
(\cite{KM1,KM2,DM3,DM5,DM7,DM6}) based on Fourier Method. But in the
case of singular potentials it may happen that the functions
$$
u_k = e^{ikx} \quad \mbox{or} \quad \sin kx, \;\;k \in \mathbb{Z},
$$
have their $L$--images outside $L^2.$ This implies, for some
singular potentials $v,$ that we have $Lf \not \in L^2$ for {\em any
smooth} (say $C^2 -$) nonzero function $f$ (see an example in
\cite{DM18}, between (1.3) and (1.4)).

In general, for any reasonable bc, the eigenfunctions $\{u_k\}$ of
the free operator $L^0_{bc}$ are not necessarily in the domain of
$L_{bc}.$  Yet, in \cite{DM17,DM16} we gave a justification of the
Fourier method for operators $L_{bc}$ with $H^{-1}$--potentials and
$bc = Per^\pm $ or $Dir.$ Our results are announced in \cite{DM17},
and in \cite{DM16} all technical details of justification of the
Fourier method are provided.

Now, in the case of singular potentials, we want to compare the
Riesz projections $P_n$ of the operator $L_{bc}, $ defined for large
enough $n$  (say $n>N$) by the formula
\begin{equation}
\label{03}
 P_n = \frac{1}{2\pi i} \int_{C_n}
(z-L_{bc})^{-1} dz, \quad C_n=\{|z-n^2|= n\},
\end{equation}
with the corresponding Riesz projections $P_n^0$ of the free
operator $L_{bc}^0$  (although $E_n^0 = Ran (P_n^0)$ maybe have no
common nonzero vectors with the domain of $L_{bc}).$ In \cite{DM18},
Theorem \ref{thm1}, we showed that
\begin{equation}
\label{03a}  \|P_n -P_n^0 \|_{L^1 \to L^\infty} \to 0 \quad
\text{as} \;n \to \infty.
\end{equation}

In this paper, the main result is Theorem \ref{thm1}, which claims,
for sufficiently large $N,$ that
\begin{equation}
\label{04} \sum_{n>N} \|P_n - P_n^0\|^2_{HS} < \infty.
\end{equation}

For a potential $v \in L^2 $ (\ref{04}) is ''easy''. Indeed, using
(\ref{011}) and (\ref{012}) below, and estimating, as in the proof
of Lemma 23 in \cite{DM15}, the Hilbert--Schmidt norm of
$VR^0_\lambda$ for $\lambda \in C_n $ (where $V$ is the operator of
multiplication by $v$ and $R^0_\lambda$ is the resolvent of the free
operator), one could get
$$
\|P_n - P_n^0\|_{HS} \leq \frac{C}{n} \|v\|_{L^2}, \quad n\geq
N(\|v\|_{L^2}),
$$
with $C$ being an absolute constant,  so (\ref{04}) follows.
However, for singular potentials $v$ the proof of (\ref{03a}) and
Theorem \ref{thm1} now is rather complicated.

Since the Hilbert--Schmidt norm does not exceed the $L^2$-norm,
(\ref{04}) implies that
\begin{equation}
\label{08}
 \sum_{n>N} \| P_n -P_n^0\|^2_{L^2 \to L^2 } <\infty,
\end{equation}
which was proven earlier by A. Savchuk and A. Shkalikov
(\cite{SS03}, Sect.2.4). This implies (by Bari--Markus theorem --
see \cite{GK}, Ch.6, Sect.5.3, Theorem 5.2) that the spectral
decompositions
\begin{equation}
\label{09} f = f_N + \sum_{n>N} P_n f
\end{equation}
converge unconditionally.

The proof of Theorem \ref{thm1} is based on the perturbation theory
(for example, see \cite{Kato}), which gives the representation
\begin{equation}
\label{011}
 P_n -P_n^0=
\frac{1}{2\pi i} \int_{C_n} \left ( R(\lambda) -R^0 (\lambda) \right
) d\lambda,
\end{equation}
where $R(\lambda) = (\lambda - L_{bc})^{-1} $ and $R^0 (\lambda )$
are the resolvents of $L_{bc}$ and of the free operator $L^0_{bc},$
respectively.

In many respects the constructions of this paper are parallel to
constructions in \cite{DM18}, the proof of (\ref{03a}); see, for
example, comments in the next paragraph. However, there is no direct
way to use the inequalities proven in \cite{DM18} and to come to the
main results of the present paper.

In the classical case,  where $v \in L^2, $ one can get reasonable
estimates for the norms $ \| R(\lambda ) - R^0 (\lambda )\|$ on the
contour $C_n,$ and then by integration for $\|P_n -P_n^0 \|.$ But
now, with $v \in H^{-1}, $ we use the same approach as in
\cite{DM18}, namely, we get good estimates for the norms $\|P_n
-P_n^0 \|$ {\em after} having integrating term by term the series
representation
\begin{equation}
\label{012} R-R^0 = R^0VR^0  +R^0VR^0 V R^0  + \cdots .
 \end{equation}
This integration kills  many terms, maybe in their matrix
representation. Only then we go to the norm estimates which allow us
to prove our main result Theorem \ref{thm1}.

\section{Main result}

 By our Theorem 21 in \cite{DM16} (about spectra localization), the operator
$L_{Per\pm}$ has, for sufficiently large $n,$ exactly two
eigenvalues (counted with their algebraic multiplicity) inside the
disc of radius $n$ about $n^2$ (periodic for even $n$ or
antiperiodic for odd $n$). The operator $L_{Dir}$ has one eigenvalue
in every such disc for all sufficiently large $n.$

Let $E_n$  be the corresponding Riesz invariant subspace, and let
$P_n$  be the corresponding Riesz projection, i.e., $$ P_n =
\frac{1}{2\pi i} \int_{C_n} (\lambda - L)^{-1} d\lambda, $$ where
$C_n = \{\lambda: \; |\lambda - n^2| =n. \}   $ Further $P_n^0$
denotes the Riesz projections of the free operator and
$\|\cdot\|_{HS}$ denotes the Hilbert--Schmidt norm.

\begin{Theorem}
\label{thm1} In the above notations, for boundary conditions $bc
=Per^\pm $ or $Dir,$
\begin{equation}
\label{110} \sum_{n>N} \|P_n - P_n^0\|^2_{HS} < \infty,
\end{equation}
\end{Theorem}

\begin{proof}
We give a complete proof in the case $bc =Per^\pm.$ If $bc =Dir$ the
proof is the same, and only minor changes are necessary due to the
fact that in this case the orthonormal system of eigenfunctions of
$L^0$ is $\{ \sqrt{2} \sin nx, \; n\in \mathbb{N}\}$ ( while it is
$\{\exp(imx), \; m\in 2 \mathbb{Z}\}$ for $bc = Per^+, $ and
$\{\exp(imx), \; m\in 1+ 2 \mathbb{Z}\}$ for $bc = Per^- $). So,
roughly speaking, the only difference is that when working with $bc
= Per^\pm $ the summation indexes in our formulas below run,
respectively, in $2 \mathbb{Z}$ and $1+ 2 \mathbb{Z},$ while for $bc
=Dir$ the summation indexes have to run in $\mathbb{N}.$ Therefore,
we consider in detail only $bc= Per^\pm. $

Now we present the proof of the theorem up to a few technical
inequalities proved in Section 3, Lemmas \ref{lemt1}, \ref{lemt2}
and \ref{lemt3}.

In \cite{DM16}, Section 5, we gave a detailed analysis of the
representation
\begin{equation}
\label{110a}  R_\lambda - R_\lambda^0 = \sum_{s=0}^\infty K_\lambda
(K_\lambda V K_\lambda )^{s+1} K_\lambda,
\end{equation}
where $K_\lambda = \sqrt{R^0_\lambda} $ -- see \cite{DM16},
(5.13-14) and what follows there.

With this definition the operator valued function $K_\lambda$ is
analytic in $ \mathbb{C} \setminus \mathbb{R}_+$. But (\ref{110a}),
(\ref{111a}) below and all formulas of this section -- which are
essentially variations of (1.11) -- have always even powers of
$K_\lambda$ and $K_\lambda^2 = R_\lambda^0 $ is analytic outside on
the complement of $ Sp(L^0).$ Certainly, this justifies the use of
Cauchy formula or theorem when warranted.

By (\ref{011}),
\begin{equation}
\label{111a} P_n - P_n^0 = \frac{1}{2\pi i} \int_{C_n}
\sum_{s=0}^\infty K_\lambda   (K_\lambda V  K_\lambda )^{s+1}
K_\lambda d\lambda
\end{equation}
if the series on the right converges. Taking into account that the
adjoint operator of $R_\lambda (v)$ is
$$
(R_\lambda (v))^* = R_{\overline{\lambda}} (\overline{v}),
$$
we get
$$ (P_n - P_n^0)^* = \frac{1}{2\pi i} \int_{C_n} \sum_{t=0}^\infty
K_\mu   (K_\mu \tilde{V}  K_\mu )^{t+1} K_\mu d\mu,
$$
where
\begin{equation}
\label{111} \tilde{V}(m) = \overline{V} (-m).
\end{equation}

 Since  $\|(P_n - P_n^0) e_m\|^2 =\langle (P_n -
P_n^0)^*(P_n - P_n^0) e_m,e_m \rangle,$ it follows that
\begin{equation}
\label{112}  \|(P_n - P_n^0) e_m\|^2 = -\frac{1}{4\pi^2 } \int
\int_{\Gamma_n} \sum_{t,s=0}^\infty \langle K_\mu (K_\mu \tilde{V}
K_\mu )^{t+1} K_\mu K_\lambda (K_\lambda V K_\lambda )^{s+1}
K_\lambda e_m, e_m \rangle  d\lambda d\mu ,
\end{equation}
where $\Gamma_n = C_n\times C_n.$ Thus,
\begin{equation}
\label{113} \sum_{n>N}\| P_n - P_n^0\|^2_{HS}=\sum_{n>N} \sum_m
 \|(P_n - P_n^0) e_m\|^2 \leq \sum_{t,s=0}^\infty A(t,s),
\end{equation}
where
\begin{equation}
\label{115}  A(t,s) = \sum_{n>N}  \left | \frac{1}{4\pi^2 } \int
\int_{\Gamma_n} \sum_m \langle K_\mu (K_\mu \tilde{V} K_\mu )^{t+1}
K_\mu K_\lambda (K_\lambda V K_\lambda )^{s+1} K_\lambda e_m, e_m
\rangle d\lambda d\mu  \right |.
\end{equation}
Notice that $A(t,s)$ depends on $N$ but this dependence is
suppressed in the notation. Our goal is to show, for sufficiently
large $N,$ that $ \sum_{t,s=0}^\infty A(t,s) <\infty $ which, in
view of (\ref{113}), implies (\ref{110}).

Let us evaluate $A(0,0).$ From the matrix representation of the
operators $K_\lambda $ and $V$ (see more details in \cite{DM16},
(5.15-22)) it follows that
\begin{equation}
\label{119} \langle K_\mu (K_\mu \tilde{V}  K_\mu )K_\mu K_\lambda
(K_\lambda V K_\lambda )K_\lambda e_m, e_m \rangle =\sum_p
\frac{\tilde{V}(m-p)V(p-m)}{(\mu -m^2)(\mu-p^2)(\lambda
-p^2)(\lambda -m^2)}.
\end{equation}
By integrating this function over $\Gamma_n = C_n \times C_n $ we
get
\begin{equation}
\label{121} \frac{1}{4\pi^2} \int \int_{\Gamma_n} \cdots =
\begin{cases} \sum_{p \neq \pm n} \frac{|V(p \mp n)|^2}{(n^2
-p^2)^2}   &  m= \pm n, \\ \frac{|V(\pm n -m)|^2}{(n^2 -m^2)^2}   &
 m \neq \pm n.
\end{cases}
\end{equation}
Thus,
$$ A(0,0)=\sum_{n>N} \sum_{p \neq \pm n} \frac{|V(p-n)|^2}{|n^2
-p^2|^2} + \sum_{n>N} \sum_{p \neq \pm n} \frac{|V(p+n)|^2}{|n^2
-p^2|^2} $$ $$+ \sum_{n>N} \sum_{m \neq \pm n}
\frac{|V(n-m)|^2}{|n^2 -m^2|^2} + \sum_{n>N} \sum_{m \neq \pm n}
\frac{|V(-n-m)|^2}{|n^2 -m^2|^2} $$ Let us estimate the first sum on
the right. In view of (\ref{001}),
\begin{equation}
\label{123} |V(m)| \leq |m| r(m), \quad  r(m)=
max(|q(m)|,|q(-m)|)\quad r \in \ell^2 (2\mathbb{Z}).
\end{equation}
Therefore, by Lemma \ref{lemt1}, we have
$$
\sum_{n>N} \sum_{p \neq \pm n} \frac{|V(p-n)|^2}{|n^2 -p^2|^2} \leq
\sum_{n>N} \sum_{p \neq \pm n} \frac{ |p-n|^2 |r(p-n)|^2}{|n^2
-p^2|^2}
$$
$$
\leq \sum_{n>N} \sum_{p \neq \pm n} \frac{|r(n-p)|^2}{|n+p|^2} \leq
C \left ( \frac{\|r\|^2}{N} + (\mathcal{E}_N (r))^2 \right ),
$$
where we use the notation
\begin{equation}
\label{125} \mathcal{E}_a (r) = \left (\sum_{|k| \geq a} |r(k)|^2
\right )^{1/2}, \quad a>0.
\end{equation}

 Since each of the other three sums could be
estimated in the same way, we get
\begin{equation}
\label{126} A(0,0) \leq C \left ( \frac{\|r\|^2}{N} + (\mathcal{E}_N
(r))^2 \right ).
\end{equation}
{\em Remark:} For convenience, here and thereafter we denote by  $C$
any absolute constant.

Next we estimate $A(t,s)$ with  $s+t >0.$ From the matrix
representation of the operators $K_\lambda $ and $V$ we get
\begin{equation}
\label{128} \langle K_\mu (K_\mu \tilde{V} K_\mu )^{t+1} K_\mu
K_\lambda (K_\lambda V K_\lambda )^{s+1} K_\lambda e_m, e_m \rangle
\end{equation}
$$
=  \sum_{i_1,\ldots,i_t,p,j_1,\ldots,j_s}
\frac{\tilde{V}(m-i_1)\tilde{V}(i_1 -i_2) \cdots \tilde{V}(i_t
-p)V(p - j_1) V(j_1-j_2) \cdots V(j_s -m)}{(\mu -m^2)(\mu -i^2_1)
\cdots (\mu -i^2_t)(\mu -p^2)(\lambda -p^2) (\lambda -j_1^2)\cdots
(\lambda -j_s^2)(\lambda -m^2)}
$$
Notice that if
\begin{equation}
\label{129}\pm n \not \in \{m,p,i_1, \ldots, i_t \} \quad  \text{or}
\quad \pm n \not \in \{m,p,j_1, \ldots, j_s \},
\end{equation}
then the integral over $C_n \times C_n $ of the corresponding term
in the above sum is zero because that term is, respectively, an
analytic function of $\mu,\; |\mu| \leq n $ and/or an analytic
function of $\lambda, \; |\lambda| \leq n.$ This observation is
crucial in finding good estimates for $A(t,s).$ It means that we may
"forget" the terms satisfying (\ref{129}).

Moreover, by the Cauchy formula, if
\begin{equation}
\label{130} m,p,i_1, \ldots, i_t \in \{\pm n\}
 \quad  \text{or} \quad
 m,p,j_1, \ldots, j_s \in \{\pm n\},
\end{equation}
then the integral of the corresponding term vanishes.

Hence we have
\begin{equation}
\label{131} A(t,s) \leq \sum_{n>N}  \left | \frac{1}{4\pi^2}
\int_{\Gamma_n} \sum_{I^*}\frac{\tilde{V}(m-i_1)\cdots \tilde{V}(i_t
-p)V(p - j_1)  \cdots V(j_s -m)}{(\mu -m^2) (\mu -i^2_1)..(\mu
-p^2)(\lambda -p^2) (\lambda -j_1^2).. (\lambda -m^2)} d\mu d\lambda
\right |,
\end{equation}
where  $I^*$ is the set of $t+s+2$-tuples of indices
$m,i_1,\ldots,i_t,p,j_1,\ldots,j_s \in n+2\mathbb{Z}$ such that
(\ref{129}) and (\ref{130}) do not hold.

In view of (\ref{131}), we may estimate $A(t,s)$ by
\begin{equation}
\label{132} A(t,s) \leq \sum_{n>N} n^2 \sup_{(\mu,\lambda) \in
\Gamma_n}\sum_{I^*} B(\mu, m,i_1,\ldots,i_t,p) \cdot  B(\lambda,
p,j_1,\ldots,j_s,m),
\end{equation}
where
\begin{equation}
\label{133} B(z, m,i_1,\ldots,i_t,p)= \frac{W(m-i_1)W(i_1 -i_2)
\cdots W(i_{t-1}-i_t) W(i_t -p)}{|z -m^2||z -i^2_1| \cdots |z
-i^2_t| |z -p^2 |},
\end{equation}
and
\begin{equation}
\label{134} B(z,m,p)= \frac{W(m-p)}{|z -m^2| |z -p^2 |}
\end{equation}
(in the degenerate case, when there are no $i$-indices), with
\begin{equation}
\label{135} W(m) = \max \{|V(m)|, |V(-m)|\}, \quad m \in
2\mathbb{Z}.
\end{equation}
In view of (\ref{123}) and (\ref{111}), we have
\begin{equation}
\label{136} W(m) = |m| r(m),  \quad \text{where} \;\; r(-m)
=r(m)\geq 0, \;\; r= (r(m)) \in 2\mathbb{Z}.
\end{equation}

We consider the following subsets of $I^*:$
\begin{equation}
\label{137}
 I_0^* = \left \{(m,i_1,\ldots,i_t,p,j_1,\ldots,j_s): \quad m=\pm n,
 \;p= \pm n \right\},
\end{equation}
\begin{equation}
\label{138}
 I_1^* = \left \{(m,i_1,\ldots,i_t,p,j_1,\ldots,j_s): \quad m=\pm n,
 \;p \neq \pm n \right\},
\end{equation}
\begin{equation}
\label{139}
 I_2^* = \left \{(m,i_1,\ldots,i_t,p,j_1,\ldots,j_s): \quad m \neq \pm n,
 \;p= \pm n \right\},
\end{equation}
\begin{equation}
\label{140}
 I_3^* = \left \{(m,i_1,\ldots,i_t,p,j_1,\ldots,j_s): \quad m \neq \pm n,
 \;p \neq \pm n \right\}.
\end{equation}
Since $I^* = \cup I^*_k, \; k=0,1,2,3,$ we have
\begin{equation}
\label{142} A(t,s) \leq A_0 (t,s)+A_1 (t,s)+A_2 (t,s)+A_3 (t,s),
\end{equation}
where $A_k (t,s)$ is the subsum of the sum on the right of
(\ref{131}) which corresponds to $I^*_k,$ i.e.,
\begin{equation}
\label{143} A_k (t,s) = \sum_{n>N} n^2 \cdot \sup_{(\mu, \lambda)
\in \Gamma_n} \sum_{I_k^*} B(\mu, m,i_1,\ldots,i_t,p) \cdot
B(\lambda, p,j_1,\ldots,j_s,m), \quad k=0,1,2,3.
\end{equation}

Let $\overline{K}_z$ denote the operator with a matrix
representation
$$(\overline{K}_z)_{jm} = \frac{1}{|z-m^2|^{1/2} \delta_{jm}}, $$
and let $W$ denote the operator with a matrix representation
$$ W_{jm} = W(j-m). $$ Then the matrix representation of the
operator $\overline{K}_z W \overline{K}_z$ is
\begin{equation}
\label{147} (\overline{K}_z W \overline{K}_z)_{jm} = \frac{W(j-m)}
{|z-j^2|^{1/2}|z-m^2|^{1/2}},
\end{equation}
and we have (see the proof of Lemma 19 in \cite{DM16})
\begin{equation}
\label{148} \|\overline{K}_z\| = \frac{1}{\sqrt{n}}, \quad
\|\overline{K}_z W \overline{K}_z\|_{HS} \leq \rho_n \quad
\text{for} \;\; z\in C_n, \; \; n \geq 3,
\end{equation}
where
\begin{equation}
\label{146} \rho_n =C \left (\mathcal{E}_{\sqrt{n}} (r) + \|r\|^2/n
\right )^{1/2},
\end{equation}
and $\|\cdot \||{HS}$ means the Hilbert--Schmidt norm of the
corresponding operator.

 Moreover, by (\ref{133}), we have
\begin{equation}
\label{149} \sum_{i_1,\ldots,i_t} B(z,m, i_1,\ldots,i_t,p) = \langle
\overline{K}_z (\overline{K}_z W \overline{K}_z)^{t+1}\overline{K}_z
e_p, e_m \rangle
\end{equation}

{\em Estimates for $A_0 (t,s).$} Notice, that $A_0 (t,0) =0$ and
$A_0 (0,s) = 0 $ because the corresponding set of indices $I^*_0$ is
empty (see the text around (\ref{130}), and the definition of
$I^*$).

Assume that $t>0, s>0.$ In view of (\ref{137}) and (\ref{143}), we
have
\begin{equation}
\label{150} A_0 (t,s) \leq \sum_{n>N} n^2 \sum_{m,p \in \{\pm n\}}
\sup_{(\mu,\lambda) \in \Gamma_n} \sum_{i_1,\ldots,i_t} B(\mu,m,
i_1,\ldots,i_t,p) \sum_{j_1,\ldots,j_s} B(\lambda,p,
j_1,\ldots,j_s,m).
\end{equation}
Therefore, by the Cauchy inequality,
\begin{equation}
\label{151} A_0 (t,s) \leq \left (\sum_{n>N} n^2 \sum_{m,p \in \{\pm
n\}} \sup_{\mu \in C_n} \left |\sum_{i_1,\ldots,i_t} B(\mu,m,
i_1,\ldots,i_t,p) \right |^2 \right )^{1/2}
\end{equation}
$$
\times \left(\sum_{n>N} n^2 \sum_{m,p \in \{\pm n\}} \sup_{\lambda
\in C_n}  \left |\sum_{j_1,\ldots,j_s} B(\lambda,p,
j_1,\ldots,j_s,m) \right |^2 \right )^{1/2}.
$$

\begin{Lemma}
\label{lem3} In the above notations,
\begin{equation}
\label{152} \sum_{n>N} n^2 \sup_{\mu \in C_n} \left |\sum_{i_1,
\ldots,i_t} B(\mu,m,i_1, \ldots,i_t,p) \right |^2 \leq C \|r\|^2
\rho_N^{2t} \quad \text{if} \quad m,p \in \{\pm n\},
\end{equation}
where $C$ is an absolute constant and $\rho_N$ is defined in
(\ref{146}).
\end{Lemma}

\begin{proof} If $t=1,$ then, by (\ref{133}),
the sum $\sigma $ in (\ref{152}) has the form
$$\sigma (m,p)= \sum_{n>N} n^2\sup_{\mu \in C_n} \left | \sum_i
\frac{W(m-i)W(i-p)}{|\mu -m^2||\mu - i^2||\mu -p^2|}  \right |^2,
\quad m,p \in \{\pm n\}.
$$
One can easily see that
$$
\sigma (-n,-n) = \sigma (n,n), \quad \sigma (-n,n) = \sigma (n,-n)
$$
by changing $i$ to $ -i$ and using that $W(-k) = W(k).$

Taking into account that $|\mu -n^2| = n $ for $\mu \in C_n,$ and
$W(k) = |k|r(k),$ we get, by the elementary inequality
\begin{equation}
\label{154} \frac{1}{|\mu - i^2|} \leq \frac{2}{|n^2 - i^2|} \quad
\text{for} \quad  \mu \in C_n, \; i \in n + 2 \mathbb{Z}, \; i \neq
\pm n,
\end{equation}
that
$$
\sigma (n,n)\leq 4\sum_{n>N} n^2 \left (\sum_{i\neq \pm n}
\frac{|n-i|}{n^2 |n +i|} r(n-i)r(i-n) + \frac{4}{n}
r(2n)r(-2n)\right )^2
$$
Therefore, by the Cauchy inequality,
$$
\sigma (n,n)\leq 4\sum_{n>N} 2n^2 \left (\sum_{i\neq \pm n}
\frac{|n-i|}{n^2 |n +i|} r(n-i)r(i-n)  \right )^2 +128
\sum_{n>N}|r(2n)r(-2n)|^2
$$
$$
\leq 2 \|r\|^2 \sum_{n>N}  \sum_{i\neq \pm n} \frac{|n-i|^2}{n^2 |n
+i|^2} |r(n-i)|^2   +128\|r\|^2 \sum_{n>N}|r(2n)|^2 \leq C \rho_N^2.
$$
(by (\ref{t3}) in Lemma \ref{lemt1}). In an analogous way, we get
$$\sigma (n,-n)= \sum_{n>N} n^2\left | \sum_{i\neq \pm n}
\frac{W(n-i)W(i+n)}{n^2|n^2 - i^2|}  \right |^2
$$
$$
= \sum_{n>N} \frac{1}{n^2} \left ( \sum_{i\neq \pm n} r(n-i)r(i+n)
\right )^2 \leq \frac{4}{N} \|r\|^4 \leq 4\rho_N^2.
$$
This completes the proof of (\ref{152}) for $t=1.$

Next we consider the case  $t>1.$
 Since $|\mu -n^2|= n$ for $\mu \in C_n,$
by (\ref{133}) the sum $\sigma $ in (\ref{152}) can be written in
the form
$$
\sigma= \sum_{n>N} \frac{1}{n^2} \sup_{\mu \in C_n} \left
|\sum_{i_1, \ldots,i_t} \frac{W(m-i_1)W(i_1-i_2)\cdots W(i_t
-p)}{|\mu-i_1^2||\mu-i_1^2|\cdots |\mu-i_t^2|}\right |^2,  \quad m,p
\in \{\pm n\}.
$$
In view of (\ref{147}), we have (with $i=i_1, k=i_t$)
$$
\sigma = \sum_{n>N} \frac{1}{n^2} \sup_{\mu \in C_n} \left
|\sum_{i,k} \frac{W(m-i)}{|\mu - i^2|^{1/2}}\cdot H_{ik} (\mu) \cdot
\frac{W(k -p)}{|\mu - k^2|^{1/2}} \right |^2, \quad m,p \in \{\pm
n\},
$$
where $ (H_{ik}(\mu))$ is the matrix representation of the operator
$ H(\mu) = (\overline{K}_\mu W \overline{K}_\mu)^{t-1}. $ By
(\ref{148}),
$$
\|H(\mu)\|_{HS} =\left ( \sum_{i,k} |H_{ik}(\mu)|^2 \right )^{1/2}
\leq \|\overline{K}_\mu W \overline{K}_\mu \|_{HS}^{t-1} \leq
\rho_N^{t-1} \quad  \text{for} \; \mu \in C_n, \; n>N. $$ Therefore,
the Cauchy inequality implies
$$
\sigma (m,p) \leq \rho_N^{2(t-1)}\cdot\sum_{n>N} \frac{1}{n^2}
\sup_{\mu \in C_n} \sum_{i,k} \frac{|W(m-i)|^2}{|\mu - i^2|}\cdot
\frac{|W(k -p)|^2}{|\mu - k^2|}.
$$
By (\ref{154}) and $W(-k) = W(k),$ one can easily see (changing $i$
with $-i, $ if necessary)  that
$$
\max_{m=\pm n} \sup_{\mu \in C_n}\sum_{i} \frac{|W(m-i)|^2}{|\mu -
i^2|} \leq \sum_{i \neq \pm n}\frac{2|W(n-i)|^2}{|n^2 - i^2|} +
\frac{|W(2n)|^2}{n}
$$
In an analogous way, it follows that
$$
\max_{p=\pm n} \sup_{\mu \in C_n}\sum_{k} \frac{|W(k-p)|^2}{|\mu -
k^2|} \leq \sum_{i \neq \pm n}\frac{2|W(n-i)|^2}{|n^2 - i^2|} +
\frac{|W(2n)|^2}{n}.
$$
Therefore, we have
$$
\sigma (m,p) \leq \rho_N^{2(t-1)}\cdot\sum_{n>N} \frac{1}{n^2} \left
(\sum_{i \neq \pm n}\frac{2|W(n-i)|^2}{|n^2 - i^2|} +
\frac{|W(2n)|^2}{n} \right )^2.
$$
Since $W(k) = |k|r(k),$  by $(a+b)^2 \leq 2a^2 +2b^2 $ and the
Cauchy inequality, we get
$$
\left (\sum_{i \neq \pm n}\frac{2|W(n-i)|^2}{|n^2 - i^2|} +
\frac{|W(2n)|^2}{n} \right )^2 \leq 8\left (  \sum_{i\neq \pm n}
\frac{|n-i|}{|n + i|} |r(n-i)|^2  \right )^2 + 32 n^2 |r(2n)|^4
$$
$$
 \leq 8\|r\|^2  \sum_{i\neq \pm n}
\frac{|n-i|^2}{|n + i|^2} |r(n-i)|^2   + 32 n^2 |r(2n)|^2 \|r\|^2.
$$
Thus,
$$
\sigma (m,p) \leq 32\|r\|^2 \rho_N^{2(t-1)} \left (\sum_{n>N}
\sum_{i\neq \pm n} \frac{|n-i|^2}{n^2 |n + i|^2} |r(n-i)|^2 +
\sum_{n>N} |r(2n)|^2 \right ) \leq C \|r\|^2 \rho_N^{2t}
$$
(by (\ref{t3}) in Lemma \ref{lemt1}).
\end{proof}

Now, by (\ref{151}) and (\ref{152}) in Lemma \ref{lem3}, we get
\begin{equation}
\label{156} A_0 (t,s) \leq C \|r\|^2 \rho_N^{t+s}, \quad t+s > 0,
\end{equation}
where $C$ is an absolute constant.

{\em Estimates for $A_1 (t,s).$} Assume that $t+ s>0.$ In view of
(\ref{138}) and (\ref{143}), we have
\begin{equation}
\label{160} A_1 (t,s) \leq \sum_{n>N} n^2 \sum_{m =\pm n,p\neq \pm
n} \sup_{\mu \in C_n} \sum_{i_1,\ldots,i_t} B(\mu,m,
i_1,\ldots,i_t,p) \sup_{\lambda \in C_n}\sum_{j_1,\ldots,j_s}
B(\lambda,p, j_1,\ldots,j_s,m).
\end{equation}
Therefore, by the Cauchy inequality,
\begin{equation}
\label{161} A_1 (t,s) \leq \left (\sum_{n>N} n^2 \sum_{m =\pm
n,p\neq \pm n} \sup_{\mu \in C_n} \left |\sum_{i_1,\ldots,i_t}
B(\mu,m, i_1,\ldots,i_t,p) \right |^2 \right )^{1/2}
\end{equation}
$$
\times \left(\sum_{n>N} n^2 \sum_{m =\pm n,p\neq \pm n}
\sup_{\lambda \in C_n}  \left |\sum_{j_1,\ldots,j_s} B(\lambda,p,
j_1,\ldots,j_s,m) \right |^2 \right )^{1/2}.
$$

\begin{Lemma}
\label{lem4} In the above notations,
\begin{equation}
\label{162} \sum_{n>N,  p \neq \pm n} n^2  \sup_{\mu \in C_n} \left
|\sum_{i_1, \ldots,i_t} B(\mu,m,i_1, \ldots,i_t,p) \right |^2 \leq C
\|r\|^2 \rho_N^{2t} \quad \text{if} \quad m \in \{\pm n\},
\end{equation}
where $C$ is an absolute constant and $\rho_N$ is defined in
(\ref{146}).
\end{Lemma}

\begin{proof} If $t=0,$ then, by (\ref{134}),
the sum $\sigma $ in (\ref{162}) has the form
$$\sigma (m)= \sum_{n>N,  p \neq \pm n} n^2\sup_{\mu \in C_n}
\frac{|W(m-p)|^2}{n^2|\mu -p^2|^2}  , \quad m = \pm n.
$$
By (\ref{154}), and since $W(-k)= W(k) =|k|r(k),$
$$\sigma (m) \leq \sum_{n>N,  p \neq \pm n}
\frac{4|W(m-p)|^2}{|n^2 -p^2|^2} = \sum_{n>N,  p \neq \pm n}
\frac{4|W(n-p)|^2}{|n^2 -p^2|^2}
$$
$$
=4\sum_{n>N,  p \neq \pm n} \frac{|r(n-p|^2}{|n+p|^2} \leq C
\rho_N^2
$$
by (\ref{t1}) in Lemma \ref{lemt1}. So, (\ref{162}) holds for $t=0.$

If $t=1,$  then, by (\ref{133}), the sum $\sigma $ in (\ref{162})
has the form
$$
\sum_{n>N,  p \neq \pm n} n^2\sup_{\mu \in C_n} \left |\sum_k
\frac{W(m-k)W(k-p)}{n|\mu -k^2||\mu -p^2|} \right |^2 , \quad m =
\pm n.
$$
By (\ref{154}), and since $W(-k)= W(k) =|k|r(k),$ we have
$$
\sigma (\pm n) \leq  \sum_{n>N,  p \neq \pm n}    \left (
\sum_{k\neq \pm n} \frac{4|n-k||k-p|}{|n^2
-k^2||n^2-p^2|}r(n-k)r(k-p) + \frac{4r(2n)r(n+p)}{|n-p|} \right )^2
$$
(to get this estimate for $m=-n$ one may replace $k$ and $p,$
respectively, by $-k$ and $-p$). Since $(a+b)^2 \leq 2a^2 +2 b^2,$
we have
$$
\sigma (\pm n) \leq 32\sigma_1 + 32\sigma_2,
$$
where
$$
\sigma_1=  \sum_{n>N,  p \neq \pm n}  \left ( \sum_{k \neq \pm n}
\frac{|k-p|}{|n+k||n^2-p^2|}r(n-k)r(k-p)  \right )^2
$$
and
$$
\sigma_2=  \sum_{n>N,  p \neq \pm n}
\frac{|r(2n)|^2|r(n+p)|^2}{|n-p|^2} \leq \|r\|^2 \cdot \sum_{n>N, p
\neq \pm n}\frac{|r(n+p)|^2}{|n-p|^2} \leq C\|r\|^2 \rho_N^2
$$
by (\ref{t1}) in Lemma \ref{lemt1}. On the other hand, the identity,
$$
\frac{k-p}{(n+k)(n+p)}= \frac{1}{n+p}-\frac{1}{n+k}
$$
implies that
$$
\sigma_1=\sum_{n>N,  p \neq \pm n}  \left ( \sum_{k \neq \pm n}
\left | \frac{1}{n+p}-\frac{1}{n+k} \right | \frac{1}{|n-p|}
r(n-k)r(k-p) \right )^2 \leq 2 \sigma_1^\prime + 2 \sigma_1^{\prime
\prime},
$$
where
$$
\sigma_1^\prime = \sum_{n>N,  p \neq \pm n} \frac{1}{|n^2-p^2|^2}
\left ( \sum_{k \neq \pm n} r(n-k)r(k-p) \right )^2 \leq \sum_{n>N,
p \neq \pm n} \frac{1}{|n^2-p^2|^2}\|r\|^2 \leq C\frac{\|r\|^2}{N},
$$
and
$$
\sigma_1^{\prime \prime}=\sum_{n>N,  p \neq \pm n}  \left ( \sum_{k
\neq \pm n} \frac{r(n-k)r(k-p)}{|n+k||n-p|}  \right )^2
$$
$$\leq
\sum_{n>N,  p \neq \pm n}  \left ( \sum_{k \neq \pm n}
\frac{|r(k-p)|^2}{|n+k|^2|n-p|^2}  \right )\cdot \|r\|^2 \leq
C\|r\|^2 \rho_N^2
$$
(by the Cauchy inequality and (\ref{t2}) in Lemma \ref{lemt1}). So,
the above inequalities imply (\ref{162}) for $t=1.$

Next we consider the case  $t>1.$
 Since $|\mu -n^2|= n$ for $\mu \in C_n,$
by (\ref{133}) the sum $\sigma $ in (\ref{162}) can be written in
the form
$$
\sigma (m)= \sum_{n>N,  p \neq \pm n}  \sup_{\mu \in C_n} \left
|\sum_{i_1, \ldots,i_t} \frac{W(m-i_1)W(i_1-i_2)\cdots W(i_t
-p)}{|\mu-i_1^2||\mu-i_1^2|\cdots |\mu-i_t^2||\mu-p^2|}\right |^2,
\quad m=\pm n.
$$
In view of (\ref{147}), we have (with $i=i_1, k=i_t$)
$$
\sigma (m) = \sum_{n>N,  p \neq \pm n}  \sup_{\mu \in C_n} \left
|\sum_{i,k} \frac{W(m-i)}{|\mu - i^2|^{1/2}}\cdot H_{ik} (\mu) \cdot
\frac{W(k -p)}{|\mu - k^2|^{1/2}|\mu-p^2|} \right |^2, \quad m=\pm
n,
$$
where $ (H_{ik}(\mu))$ is the matrix representation of the operator
$ H(\mu) = (\overline{K}_\mu W \overline{K}_\mu)^{t-1}. $ By
(\ref{148}),
$$
\|H(\mu)\|_{HS} =\left ( \sum_{i,k} |H_{ik}(\mu)|^2 \right )^{1/2}
\leq \|\overline{K}_\mu W \overline{K}_\mu \|_{HS}^{t-1} \leq
\rho_N^{t-1} \quad  \text{for} \; \mu \in C_n, \; n>N. $$ Therefore,
the Cauchy inequality and (\ref{154}) imply
$$
\sigma (\pm n) \leq 4\rho_N^{2(t-1)}\cdot\sum_{n>N,  p \neq \pm n}
\frac{1}{ (n^2-p^2)^2} \sup_{\mu \in C_n} \sum_{i,k}
\frac{|W(n+i)|^2}{|\mu - i^2|}\cdot \frac{|W(k + p)|^2}{|\mu - k^2|}
$$
(one may see that the inequality holds for $m=\pm n$ by replacing,
if necessary, $i $ by $-i$ and $p$ by $-p $).

From (\ref{154}) and $W(k)= |k|r(k)$ it follows that
$$
\sup_{\mu \in C_n} \sum_{i} \frac{|W(n+i)|^2}{|\mu - i^2|}  \leq
2\sum_{i\neq \pm n} \frac{|n+i|}{|n-i|}r(n+i)|^2  + 4n |r(2n)|^2
$$
and
$$
\sup_{\mu \in C_n} \sum_{k} \frac{|W(k + p)|^2}{|\mu - k^2|} \leq
2\sum_{k \neq \pm n} \frac{|k+p|^2}{|n^2-k^2|}|r(k+p|^2 +
\frac{|n+p|^2}{n} |r(n+p|^2 +\frac{|n-p|^2}{n} |r(n-p|^2.
$$
Therefore, we have
$$
\sigma (\pm n) \leq 4 \rho_N^{2(t-1)} ( 4\sigma_1+2\sigma_2
+2\sigma_3 + 8\sigma_4 +4\sigma_5 +4\sigma_6),
$$
where
$$
\sigma_1 =
 \sum_{n>N, p\neq \pm n} \frac{1}{|n^2-p^2|^2} \sum_{i,k
\neq \pm n} \frac{|n+i||p+k|^2}{|n-i||n^2-k^2|} |r(n+i)|^2|r(p+k)|^2
 \leq C \|r\|^2 \rho_N^2
$$
(by Lemma \ref{lemt2});
$$
\sigma_2 =
 \sum_{n>N, p\neq \pm n} \frac{|n+p|^2}{|n^2-p^2|^2}|r(n+p)|^2 \sum_{i
\neq \pm n} \frac{|n+i|}{n|n-i|} |r(n+i)|^2
 $$
$$
\leq  \sum_{n>N, p\neq \pm n} \frac{|r(n+p)|^2}{|n-p|^2} \cdot
2\|r\|^2 \leq C \|r\|^2 \rho_N^2
$$
(since $\frac{|n+i|}{n|n-i|} = \left |\frac{1}{n-i} -\frac{1}{2n}
\right | \leq 2,$ and by (\ref{t1}) in Lemma \ref{lemt1});
$$
\sigma_3 =
 \sum_{n>N, p\neq \pm n} \frac{|n-p|^2}{|n^2-p^2|^2}|r(n-p)|^2 \sum_{i
\neq \pm n} \frac{|n+i|}{n|n-i|} |r(n+i)|^2 =\sigma_2 \leq  C
\|r\|^2 \rho_N^2
 $$
( the change $p \to -p $ shows that $\sigma_3 = \sigma_2 $);
$$
\sigma_4 =
 \sum_{n>N, p\neq \pm n} \frac{n}{|n^2-p^2|^2}|r(2n)|^2
\sum_{k \neq \pm n} \frac{|k+p|^2}{|n^2-k^2|}|r(k+p|^2 \leq  C
\|r\|^2 \rho_N^2
$$
(by Lemma \ref{lemt3};
$$
\sigma_5 =
 \sum_{n>N, p\neq \pm n} \frac{|n+p|^2}{|n^2-p^2|^2}|r(2n)|^2|r(n+p)|^2
\leq \sum_{n>N} |r(2n)|^2 \sum_{ p\neq \pm n}|r(n+p)|^2\leq  C
\|r\|^2 \rho_N^2
$$
and
$$
\sigma_6 =
 \sum_{n>N, p\neq \pm n} \frac{|n-p|^2}{|n^2-p^2|^2}|r(2n)|^2|r(n-p)|^2
= \sigma_5 \leq  C \|r\|^2 \rho_N^2
$$
(the change $p \to -p $ shows that $\sigma_6 = \sigma_5 $). Hence
$$
\sigma (\pm n) \leq C \|r\|^2 \rho_N^{2t},
$$
which completes the proof of (\ref{162}).
\end{proof}

Now, by (\ref{161}) and (\ref{162}) in Lemma \ref{lem4}, we get
\begin{equation}
\label{166} A_1 (t,s) \leq C \|r\|^2 \rho_N^{t+s}, \quad t+s > 0,
\end{equation}
where $C$ is an absolute constant.

{\em Estimates for $A_2 (t,s).$} Since $m$ and $p$ play  symmetric
roles, the same argument that was used to estimate $A_1 (t,s)$
yields
\begin{equation}
\label{168} A_2 (t,s) \leq C \|r\|^2 \rho_N^{t+s}, \quad t+s > 0,
\end{equation}
where $C$ is an absolute constant.

{\em Estimates for $A_3 (t,s).$}  In view of (\ref{140}) and the
definition of the set $I^*$ (see the text after (\ref{131})),
$I^*_3$ is the set of $t+s+2$-tuples of indices
$(m,i_1,\ldots,i_t,p,j_1,\ldots,j_s)$ such that $t \geq 1,$ $s\geq
1, $ and
$$
m,p \neq \pm n,\quad \{i_1,\ldots,i_t\}\cap \{\pm n\} \neq
\emptyset,\quad \{j_1,\ldots,j_s\}\cap \{\pm n\} \neq \emptyset.
$$
Therefore, by (\ref{143}), we have
\begin{equation}
\label{170} A_3 (t,s) \leq \sum_{n>N} n^2 \sum_{m,p \neq \pm n}
\sup_{\mu \in C_n} \sum^*_{i_1,\ldots,i_t} B(\mu,m,
i_1,\ldots,i_t,p) \sup_{\lambda \in C_n}\sum^*_{j_1,\ldots,j_s}
B(\lambda,p, j_1,\ldots,j_s,m),
\end{equation}
where $*$ means that at least one of the summation indices is equal
to $\pm n.$ The Cauchy inequality implies
\begin{equation}
\label{171} A_3 (t,s) \leq \left (\sum_{n>N} n^2 \sum_{m,p\neq \pm
n} \sup_{\mu \in C_n} \left |\sum_{i_1,\ldots,i_t} B(\mu,m,
i_1,\ldots,i_t,p) \right |^2 \right )^{1/2}
\end{equation}
$$
\times \left(\sum^*_{n>N} n^2 \sum_{m ,p\neq \pm n} \sup_{\lambda
\in C_n}  \left |\sum^*_{j_1,\ldots,j_s} B(\lambda,p,
j_1,\ldots,j_s,m) \right |^2 \right )^{1/2}.
$$

\begin{Lemma}
\label{lem5} In the above notations,
\begin{equation}
\label{172} \sum_{n>N}\sum_{m,p \neq \pm n} n^2  \sup_{\mu \in C_n}
\left |\sum^*_{i_1, \ldots,i_t} B(\mu,m,i_1, \ldots,i_t,p) \right
|^2 \leq C t \|r\|^4 \rho_N^{2(t-1)},
\end{equation}
where $C$ is an absolute constant and $\rho_N$ is defined in
(\ref{146}).
\end{Lemma}

\begin{proof}
Let $\tau \leq t$ be the least integer such that $i_\tau =\pm n.$
Then, by (\ref{133}) or (\ref{134}), and since $|\mu -n^2|=n$ for
$\mu \in C_n,$
$$
B(\mu,m,i_1,\ldots,i_{\tau -1},\pm n,i_{\tau +1}, \ldots,i_t,p)= n
B(\mu,m,i_1, \ldots,i_{\tau -1},\pm n)\cdot B(\mu,\pm n,i_{\tau +1},
\ldots,i_t,p).
$$
Therefore, if $\sigma $ denotes the sum in (\ref{172}), we have
$$
\sigma  \leq \sum_{\tau=1}^t  \sum_{\tilde{n}=\pm n}\sum_{n>N} n^4
\sum_{m \neq \pm n}\sup_{\mu \in C_n} \left
|\sum_{i_1,\ldots,i_{\tau -1}} B(\mu,m,i_1, \ldots,i_{\tau
-1},\tilde{n}) \right |^2
$$
$$ \times
 \sum_{p \neq \pm n} \sup_{\mu \in C_n} \left |
 \sum_{i_{\tau +1},
\ldots,i_t} B(\mu,\tilde{n},i_{\tau +1}, \ldots,i_t,p) \right |^2
$$

On the other hand, by Lemma \ref{lem4},
$$
n^2 \cdot \sum_{p \neq \pm n} \sup_{\mu \in C_n} \left |
 \sum_{i_{\tau +1},
\ldots,i_t} B(\mu,\tilde{n},i_{\tau +1}, \ldots,i_t,p) \right |^2
\leq  C \|r\|^2 \rho_N^{2(t-\tau)},\quad n>N.
$$
Thus, we have
$$
\sigma \leq  C \|r\|^2 \sum_{\tau=1}^t \rho_N^{2(t-\tau)}
\sum_{\tilde{n}=\pm n} \sum_{n>N} n^2 \sum_{m \neq \pm n}\sup_{\mu
\in C_n} \left |\sum_{i_1,\ldots,i_{\tau -1}} B(\mu,m,i_1,
\ldots,i_{\tau -1},\tilde{n}) \right |^2
$$
Again by Lemma \ref{lem4},
$$
\sum_{n>N} n^2 \sum_{m \neq \pm n}\sup_{\mu \in C_n} \left
|\sum_{i_1,\ldots,i_{\tau -1}} B(\mu,m,i_1, \ldots,i_{\tau
-1},\tilde{n}) \right |^2 \leq  C \|r\|^2 \rho_N^{2(\tau -1)}
$$
(one may apply Lemma \ref{lem4} because $ B(\mu,m,i_1,
\ldots,i_{\tau -1},\tilde{n}) =B(\mu,\tilde{n},j_1, \ldots,j_{\tau
-1},m)$) if $j_1 =i_{\tau -1}, \ldots, j_{\tau -1}= i_1 $). Hence,
$$
\sigma \leq  C \|r\|^4 \sum_{\tau=1}^t \rho_N^{2(t-1)}=C t \|r\|^4
 \rho_N^{2(t-1)},
$$
which completes the proof.
\end{proof}

By (\ref{171}) and (\ref{172}) (since the roles of $m$ and $p$ are
symmetric in (\ref{171})), we get
\begin{equation}
\label{180}
 A_3 (t,s) \leq C \sqrt{ts} \|r\|^4  \rho_N^{(t+s-2)}
\leq C (t+s) \|r\|^4  \rho_N^{(t+s-2)}.
\end{equation}
Now we are ready to complete the proof of Theorem \ref{thm1}. Choose
$N$ so large that $ \rho_N <1. $ Then, from (\ref{126}),
(\ref{142}), (\ref{156}), (\ref{166}), (\ref{168}) and (\ref{180})
it follows that
$$
\sum_{t,s=0}^\infty A(t,s) <\infty,
$$
which, in view of (\ref{113}), yields (\ref{110}).
\end{proof}
So, Theorem \ref{thm1} is proven subject to Lemmas
\ref{lemt1},\ref{lemt2} and \ref{lemt3} in the next section.

\section{Technical Lemmas}
Throughout this section we use that
\begin{equation}
\label{t0} \sum_{n>N}\frac{1}{n^2} < \sum_{n>N}\left (
\frac{1}{n-1}- \frac{1}{n}  \right ) = \frac{1}{N}, \quad N \geq 1.
\end{equation}
and
\begin{equation}
\label{t00} \sum_{p \neq \pm n} \frac{1}{(n^2-p^2)^2} <
\frac{4}{n^2}, \quad   n\geq 1
\end{equation}
(since
$$
\frac{1}{(n^2-p^2)^2} = \frac{1}{4n^2} \left (  \frac{1}{n-p}+
\frac{1}{n+p} \right )^2 \leq \frac{1}{2n^2} \left (
\frac{1}{(n-p)^2}+ \frac{1}{(n+p)^2} \right ),
$$
 the sum in (\ref{t00}) does not exceed
$$\frac{1}{2n^2}
\left ( \sum_{p \neq \pm n} \frac{1}{(n-p)^2}+\sum_{p \neq \pm n}
\frac{1}{(n+p)^2} \right ) \leq \frac{1}{2n^2}\cdot 2
\frac{\pi^2}{3} <\frac{4}{n^2}
$$
because $\pi^2 <10$).

\begin{Lemma}
\label{lemt1} If $r = (r(k)) \in \ell^2 (2\mathbb{Z}) $ (or $r =
(r(k)) \in \ell^2 (\mathbb{Z}) $), then
\begin{equation}
\label{t1} \sum_{n>N,k\neq n} \frac{|r(n+k)|^2}{|n-k|^2} \leq C
\left ( \frac{\|r\|^2}{N} + (\mathcal{E}_N (r))^2 \right ),
\end{equation}
\begin{equation}
\label{t2} \sum_{n>N, k\neq n} \frac{|n+k|^2}{n^2|n-k|^2}
|r(n+k)|^2\leq C \left ( \frac{\|r\|^2}{N} + (\mathcal{E}_N (r))^2
\right ),
\end{equation}
and
\begin{equation}
\label{t3} \sum_{n>N,p, k\neq n} \frac{|r(p+k)|^2}{|n-p|^2|n-k|^2}
\leq C \left ( \frac{\|r\|^2}{N} + (\mathcal{E}_N (r))^2 \right ),
\end{equation}
where $n \in \mathbb{N},\; k,p \in n+ 2\mathbb{Z} $ (or,
respectively, $k,p \in \mathbb{Z} $) and  $C $ is an absolute
constant.
\end{Lemma}

\begin{proof} Indeed, we have (with $\tilde{k} = n+k, $ and using (\ref{t0}))
$$
\sum_{n>N, k \neq n} \frac{|r(n+k)|^2}{|n-k|^2} = \sum_{n>N,k<0}
\frac{|r(n+k)|^2}{|n-k|^2} + \sum_{n>N} \sum_{0 \leq k \neq n}
\frac{|r(n+k)|^2}{|n-k|^2}
$$
$$
\leq \sum_{n>N} \frac{1}{n^2} \sum_{\tilde{k}} |r(\tilde{k})|^2 +
\sum_{\tilde{k}> N} |r(\tilde{k})|^2 \sum_{n \neq \tilde{k}/2}
\frac{1}{|2n-\tilde{k}|^2} \leq  C \left ( \frac{\|r\|^2}{N} +
(\mathcal{E}_N (r))^2 \right ).$$

Next we prove (\ref{t2}). By the identity
$$
\frac{n+k}{n(n-k)} = \frac{1}{n-k}-\frac{1}{2n},
$$
we get  (using the inequality $ab \leq (a^2 +b^2)/2$)
$$
\sum_{n>N, k\neq n} \frac{|n+k|^2}{n^2|n-k|^2} |r(n+k)|^2=
\sum_{n>N, k\neq n} \left ( \frac{1}{n-k}-\frac{1}{2n} \right )^2
|r(n+k)|^2
$$
$$
\leq \frac{1}{2} \sum_{n>N, k\neq n} \frac{|r(n+k)|^2}{|n-k|^2} +
\frac{1}{2}   \sum_{n>N} \frac{1}{4 n^2} \sum_k |r(n+k)|^2.
$$
In view of (\ref{t0}) and (\ref{t1}), from here (\ref{t2}) follows.

In order to prove (\ref{t3}), we set $\tilde{p}= n-p$ and
$\tilde{k}= n-k. $ Then
$$
\sum_{n>N; p, k\neq n} \frac{|r(p+k)|^2}{|n-p|^2|n-k|^2}=
\sum_{\tilde{p},\tilde{k} \neq 0} \frac{1}{\tilde{p}^2}
\frac{1}{\tilde{k}^2}\sum_{n>N} |r(2n-\tilde{p}-\tilde{k}|^2
$$
$$
\leq \sum_{0<|\tilde{p}|,|\tilde{k}|\leq N/2} \frac{1}{\tilde{p}^2}
\frac{1}{\tilde{k}^2} \sum_{n>N} |r(2n-\tilde{p}-\tilde{k}|^2
+\sum_{|\tilde{p}| > N/2} \sum_{|\tilde{k}|\neq 0} \cdots +
\sum_{|\tilde{p}|\neq 0} \sum_{|\tilde{k}|> N/2} \cdots
$$
$$
\leq C(\mathcal{E}_N (r))^2 + \frac{C}{N}\|r\|^2
+\frac{C}{N}\|r\|^2,
$$
which completes the proof.
\end{proof}

\begin{Lemma}
\label{lemt2} Suppose that $r = (r(k)) \in \ell^2 (2\mathbb{Z}) $
(or $r = (r(k)) \in \ell^2 (\mathbb{Z}). $) Then
\begin{equation}
\label{t9} \sum_{n>N,p\neq \pm n} \frac{1}{|n^2-p^2|^2} \sum_{i,k
\neq \pm n} \frac{|n+i||k+p|^2}{|n-i||n^2-k^2|} |r(n+i)|^2
|r(k+p)|^2
 \leq C \|r\|^2   \left ( \frac{\|r\|^2}{N} +
(\mathcal{E}_N (r))^2 \right ),
\end{equation}
where $C$ is an absolute constant.
\end{Lemma}

\begin{proof} Let $\Sigma $
be the sum in (\ref{t9}). Taking into account that
$$
\frac{k+p}{(n-p)(n+k)} =\frac{1}{n-p} -\frac{1}{n+k},\quad
\frac{k+p}{(n+p)(n-k)}=\frac{1}{n-k} -\frac{1}{n+p}
$$
and $ (n+i)/(n-i)= 2n/(n-i) -1, $ we get
$$
\Sigma \leq \sum \frac{1}{|n^2-p^2|} \left | \frac{1}{n-p}
-\frac{1}{n+k}  \right | \left | \frac{1}{n-k} -\frac{1}{n+p} \right
| \left | \frac{2n}{n-i} -1 \right ||r(n+i)|^2|r(p+k)|^2
$$
Therefore,
\begin{equation}
\label{t10} \Sigma \leq \sum_{\nu =1}^8 \Sigma_\nu,
\end{equation}
with
\begin{equation}
\label{t11} \Sigma_1 = \sum \frac{1}{|n^2-p^2|^2} \frac{2n}{|n-i|}
|r(n+i)|^2|r(k+p)|^2,
\end{equation}
\begin{equation}
\label{t12} \Sigma_2 = \sum \frac{1}{|n^2-p^2|} \frac{1}{|n^2-k^2|}
\frac{2n}{|n-i|} |r(n+i)|^2|r(k+p)|^2,
\end{equation}
\begin{equation}
\label{t13} \Sigma_3 = \sum \frac{1}{|n^2-p^2|} \frac{1}{|n-p|}
\frac{1}{|n-k|} \frac{2n}{|n-i|} |r(n+i)|^2|r(k+p)|^2,
\end{equation}
\begin{equation}
\label{t14} \Sigma_4 = \sum \frac{1}{|n^2-p^2|} \frac{1}{|n+p|}
\frac{1}{|n+k|} \frac{2n}{|n-i|} |r(n+i)|^2|r(k+p)|^2,
\end{equation}
\begin{equation}
\label{t15} \Sigma_5 = \sum \frac{1}{|n^2-p^2|^2}
|r(n+i)|^2|r(k+p)|^2,
\end{equation}
\begin{equation}
\label{t16} \Sigma_6 = \sum \frac{1}{|n^2-p^2|} \frac{1}{|n^2-k^2|}
 |r(n+i)|^2|r(k+p)|^2,
\end{equation}
\begin{equation}
\label{t17} \Sigma_7 = \sum \frac{1}{|n^2-p^2|} \frac{1}{|n-p|}
\frac{1}{|n-k|}  |r(n+i)|^2|r(k+p)|^2,
\end{equation}
\begin{equation}
\label{t18} \Sigma_8 = \sum \frac{1}{|n^2-p^2|} \frac{1}{|n+p|}
\frac{1}{|n+k|} |r(n+i)|^2|r(k+p)|^2,
\end{equation}
where the summation is over  $n>N $ and $ i,k,p \neq \pm n.$

After summation over $k$ in (\ref{t11}) we get, in view of
(\ref{t00}),
$$
\Sigma_1 \leq \|r\|^2  \cdot \sum_{n>N, i \neq \pm n}
\frac{2n}{|n-i|} |r(n+i)|^2 \sum_{p\neq \pm n} \frac{1}{|n^2-p^2|^2}
$$
$$
\leq C \|r\|^2  \cdot \sum_{n>N, i \neq \pm n} \frac{1}{|n-i|}
\frac{1}{n} |r(n+i)|^2
$$
$$
\leq C \|r\|^2  \cdot \left ( \sum_{n>N, i \neq \pm n}
\frac{|r(n+i)|^2}{|n-i|^2} + \sum_{n>N, i \neq \pm n}
\frac{|r(n+i)|^2}{n^2} \right ).
$$
From here it follows, in view of (\ref{t0}) and (\ref{t1}), that
\begin{equation}
\label{t21} \Sigma_1 \leq C_1 \|r\|^2 \left ( \frac{\|r\|^2}{N} +
(\mathcal{E}_N (r))^2 \right ).
\end{equation}
By the inequality $2ab \leq a^2 + b^2,$ considered with $a= 1/|n^2
-p^2|$ and $b= 1/|n^2 -k^2|,$ one can easily see that
\begin{equation}
\label{t22} \Sigma_2 \leq \Sigma_1.
\end{equation}
Since
$$
\frac{2n}{n^2-p^2} = \frac{1}{n-p} + \frac{1}{n+p},
$$
we have
$$\Sigma_3 \leq \Sigma_3^\prime +\Sigma^{\prime \prime}_3,
$$
where
$$\Sigma_3^\prime =
\sum  \frac{1}{|n-p|^2} \frac{1}{|n-k|} \frac{1}{|n-i|}
|r(n+i)|^2|r(k+p)|^2.
$$
and
$$\Sigma^{\prime \prime}_3 =
\sum  \frac{1}{|n^2-p^2|} \frac{1}{|n-k|} \frac{1}{|n-i|}
|r(n+i)|^2|r(k+p)|^2.
$$

 The inequality $2ab \leq a^2 + b^2,$ considered with $a= 1/|n
-k|$ and $b= 1/|n -i|,$  yields
$$\Sigma_3^\prime
\leq \frac{1}{2} \sum_{n>N;p, k\neq n}
\frac{|r(k+p)|^2}{|n-p|^2|n-k|^2} \sum_i |r(n+i)|^2 $$ $$ +
\frac{1}{2} \sum_{n>N;i\neq n} \frac{|r(n+i)|^2}{|n-i|^2}
\sum_{p\neq n} \frac{1}{|n-p|^2} \sum_k |r(k+p)|^2\leq  C \left (
\frac{\|r\|^2}{N} + (\mathcal{E}_N (r))^2 \right ) \|r\|^2
$$
(by (\ref{t1}) and (\ref{t3}) in Lemma \ref{lemt1}). In an analogous
way, by the Cauchy inequality and (\ref{t1}) and (\ref{t3}) in Lemma
\ref{lemt1}, we get
$$\Sigma^{\prime \prime}_3 \leq \left ( \sum_{n>N;p, k\neq n}
\frac{|r(k+p)|^2}{|n-p|^2|n-k|^2} \sum_i |r(n+i)|^2 \right )^{1/2}
$$
$$
\times \left ( \sum_{n>N;i\neq n} \frac{|r(n+i)|^2}{|n-i|^2}
\sum_{p\neq n} \frac{1}{|n+p|^2} \sum_k |r(k+p)|^2 \right )^{1/2}
\leq  C \left ( \frac{\|r\|^2}{N} + (\mathcal{E}_N (r))^2 \right )
\|r\|^2.
$$
Thus,
\begin{equation}
\label{t23} \Sigma_3 \leq C \left ( \frac{\|r\|^2}{N} +
(\mathcal{E}_N (r))^2 \right ) \|r\|^2.
\end{equation}

Next we estimate $\Sigma_7.$ After summation over $i $ we get
$$
\Sigma_7 = \|r\|^2 \cdot \sum \frac{1}{|n^2-p^2|} \frac{1}{|n-p|}
\frac{1}{|n-k|}  |r(p+k)|^2.
 $$
Now the Cauchy inequality implies
$$
\Sigma_7 \leq
 \|r\|^2
 \left ( \sum_{n>N;p\neq \pm n} \frac{1}{|n^2-p^2|^2}
\sum_k |r(p+k)|^2 \right )^{1/2}
 \left ( \sum_{n>N;p, k\neq n} \frac{|r(p+k)|^2}{|n-k|^2|n-p|^2}
   \right )^{1/2}
$$
Therefore, by (\ref{t0}), (\ref{t00}), and (\ref{t3}) in Lemma
\ref{lemt1},
\begin{equation}
\label{t24} \Sigma_7 \leq C \left ( \frac{\|r\|^2}{N} +
(\mathcal{E}_N (r))^2 \right ) \|r\|^2.
\end{equation}
To estimate $\Sigma_4$ and $\Sigma_8,$ notice that if $
|r(-k)|=|r(k)| \; \forall \,k $ (which we can always assume because
otherwise one may replace $(r(k))$ by $(|r(k)|+|r(-k)|)$),
 then the change of indices $p
\to -p $ and $k \to -k $ leads to $ \Sigma_4= \Sigma_3 $ and $
\Sigma_8= \Sigma_7. $ Thus
\begin{equation}
\label{t26} \Sigma_4 \leq C \left ( \frac{\|r\|^2}{N} +
(\mathcal{E}_N (r))^2 \right ) \|r\|^2, \quad \Sigma_8 \leq C \left
( \frac{\|r\|^2}{N} + (\mathcal{E}_N (r))^2 \right ) \|r\|^2.
\end{equation}
By the inequality $2ab \leq a^2 + b^2,$ considered with $a= 1/|n^2
-p^2|$ and $b= 1/|n^2 -k^2|,$ one can easily see that
\begin{equation}
\label{t30} \Sigma_6 \leq \Sigma_5.
\end{equation}
Finally, by (\ref{t0}) and (\ref{t00}), we get
\begin{equation}
\label{t31} \Sigma_5 =  \sum_{n>N;p\neq \pm n} \frac{1}{|n^2-p^2|^2}
\sum_k |r(k+p)|^2 \sum_i |r(n+i)|^2 \leq \frac{C}{N}\|r\|^4.
\end{equation}
Now, (\ref{t10})--(\ref{t31}) imply (\ref{t9}), which completes the
proof.
\end{proof}

\begin{Lemma}
\label{lemt3} In the above notations, we have
\begin{equation}
\label{t33}  \sum_{n>N, p\neq \pm n} \frac{n}{|n^2-p^2|^2}|r(2n)|^2
\sum_{k \neq \pm n} \frac{|k+p|^2}{|n^2-k^2|}|r(k+p)|^2 \leq  C
\|r\|^2 (\mathcal{E}_N (r))^2 .
\end{equation}
\end{Lemma}

\begin{proof}
Let $\Sigma $ be the sum in (\ref{t33}). The identities
$$
\frac{k+p}{(n-p)(n+k)} = \frac{1}{n-p} - \frac{1}{n+k}, \quad
\frac{k+p}{(n+p)(n-k)} = \frac{1}{n-k} - \frac{1}{n+p},
$$
and the inequality $n \leq |n^2-p^2|, \; p \neq \pm n, $ imply that
$$
\Sigma \leq \sum_{n>N} \sum_{k,p\neq \pm n} \left | \frac{1}{n-p} -
\frac{1}{n+k} \right | \left | \frac{1}{n-k} - \frac{1}{n+p} \right
| |r(2n)|^2 |r(k+p)|^2 \leq \Sigma_1 +\Sigma_2 + \Sigma_3 +\Sigma_4,
$$
where
$$
\Sigma_1 = \sum_{n>N} |r(2n)|^2 \sum_{p\neq \pm n}
\frac{1}{|n^2-p^2|} \sum_{k\neq \pm n} |r(k+p)|^2 \leq C
(\mathcal{E}_N (r))^2 \|r\|^2;
$$
$$
\Sigma_2 = \sum_{n>N} |r(2n)|^2 \sum_{k\neq \pm n}
\frac{1}{|n^2-k^2|} \sum_{p\neq \pm n} |r(k+p)|^2 \leq C
(\mathcal{E}_N (r))^2 \|r\|^2;
$$
$$
\Sigma_3=\sum_{n>N} \sum_{k,p\neq \pm n}
\frac{1}{|n-p|}\frac{1}{|n-k|}|r(2n|^2|r(k+p)|^2
$$
and
$$
\Sigma_4=\sum_{n>N} \sum_{k,p\neq \pm n}
\frac{1}{|n+p|}\frac{1}{|n+k|}|r(2n)|^2|r(k+p)|^2.
$$
The inequality $2ab \leq a^2 + b^2$ yields $\Sigma_3\leq
\Sigma_3^\prime + \Sigma_3^{\prime \prime}$ with
$$
\Sigma_3^\prime =  \sum_{n>N}|r(2n|^2 \sum_{p\neq \pm n}
\frac{1}{|n-p|^2}\sum_{k\neq \pm n}|r(k+p)|^2\leq C (\mathcal{E}_N
(r))^2 \|r\|^2
$$
and
$$
\Sigma_3^{\prime \prime}=\sum_{n>N}|r(2n|^2 \sum_{k\neq \pm n}
\frac{1}{|n-k|^2}\sum_{p\neq \pm n}|r(k+p)|^2\leq C (\mathcal{E}_N
(r))^2 \|r\|^2.
$$
Therefore,
$$
\Sigma_3 \leq C \|r\|^2 (\mathcal{E}_N (r))^2.
$$
The same argument shows that
$$
\Sigma_4 \leq C \|r\|^2 (\mathcal{E}_N (r))^2,
$$
which completes the proof.
 \end{proof}

\section{Unconditional Convergence of Spectral Decompositions}

 1. To be accurate we should mention that in Formula (\ref{09}) the first
vector-term $f_N$ is defined as $P^N f,$  where (see \cite{DM18},
(5.40))
\begin{equation}
\label{c1} P^N = \frac{1}{2\pi i} \int_{\partial R_N}
(z-L_{bc})^{-1} dz,
\end{equation}
and $R_N$ is the rectangle
\begin{equation}
\label{c2} R_N = \{ z\in \mathbb{C}: \; -N < Re z < N^2 +N, \; |Im z
| < N \}.
\end{equation}

The Bari--Markus Theorem (\cite{B,M}; \cite{GK}, Section 5.2) gives
us the claim (\ref{09}) if the following hypotheses hold:

$$ (a) \qquad \qquad
 \sum_{n>N} \| P_n -P_n^0\|^2_{L^2 \to L^2 } <\infty \quad \text{for some} \;N,$$

$$ (b)\qquad   \text{Codim} \, H_m = \text{Codim} \, H_m^0  \quad
\text{for sufficiently large} \;m,$$
   where
   $$H_m = \text{Lin Span} \{Ran P_k, \; k\geq m\}, \qquad
H_m^0 = \text{Lin Span} \{Ran P^0_k \; k\geq m\} $$

Theorem 1 implies (a). On the other hand (b) is proven in details in
\cite{DM18}, see Theorem 21, in particular, (5.54) and (5.56).
Therefore we come to the following.

\begin{Proposition}
\label{prop2} Under the conditions of Theorem \ref{thm1}, if $N$ is
sufficiently large,  then for any $f \in L^2 (I)$
\begin{equation}
\label{c09} f = P^N f + \sum_{n>N} P_n f;
\end{equation}
these series converge unconditionally  in $ L^2 (I) .$
\end{Proposition}

  This statement has been given in \cite{SS03}, Section 2.4. Our alternative
proof is based on Fourier method which has been justified in the
analysis of Hill operators with $H^{-1}$ potentials in our paper
\cite{DM16} (see \cite{DM17} as well). \vspace{3mm}

    2. In this context it is worth to mention a version of the Bari--Markus
theorem  in the case of 1D periodic Dirac operators
$$ Ly = i \begin{pmatrix} 1
& 0 \\ 0 & -1
\end{pmatrix}
\frac{dy}{dx}  + V(x) y, \quad y = \begin{pmatrix} y_1\\y_2
\end{pmatrix}, $$
where $$V(x) = \begin{pmatrix} 0 & P(x) \\ Q(x) & 0 \end{pmatrix},
\quad V(x+ \pi) = V(x), \quad P,Q \in L^2 ([0,\pi]). $$

For Riesz projections (in the case of $bc = Per^\pm$ and $ Dir$ -
see definitions and details in \cite{M04} or \cite{DM15}, Sect. 1.1)
Theorem 8.8 in \cite{M04} or Theorem 4 in \cite{M03} claims the
following:
\begin{Proposition}
\label{prop3} Let $\Omega = (\Omega (k)), \; k \in \mathbb{Z},$ be a
weight such that
\begin{equation}
\label{c11} \sum \frac{1}{(\Omega(k))^2} <\infty.
\end{equation}
If $V \in H(\Omega),$ then
\begin{equation}
\label{c12} F = P^N F + \sum_{|n|>N} P_n F    \quad \forall F \in
L^2;
\end{equation}
these spectral decompositions converge unconditionally.
\end{Proposition}

\end{document}